\documentclass[12pt,twoside]{article}

\usepackage{amssymb,amsmath}
\usepackage{latexsym}
\usepackage{changebar}

\def\Sym{{\bf Sym}}

\def\<{\langle}
\def\>{\rangle}

\def\Z{{\mathbb Z}}

\def\tHH{{\rm \tilde H}}

\def\tGG{{\rm \tilde G}}

\def\tPP{{\rm \tilde P}}

\def\Des{{\rm Des\,}}
\def\bkt{{\bf \tilde k}}
\def\K{{\bf K}}

\def\vt{\tilde{v}}
\def\ttH{{\bf \widetilde{H}}}
\def\tHH{{\rm \tilde H}}

\newdimen\squaresize
\newdimen\thickness        
                                                    
\def\square#1{\hbox{\vrule width \thickness
   \vbox to \squaresize{\hrule height \thickness\vss                            
      \hbox to \squaresize{\hss#1\hss}
   \vss\hrule height\thickness} 
\unskip\vrule width \thickness} 
\kern-\thickness}                                                            
                               
\def\vsquare#1{\vbox{\square{$\casestyle#1$}}\kern-\thickness}
\def\blank{\omit\hskip\squaresize}

\def\young#1{\vcenter{%
    \squaresize=18pt\thickness=0.5pt\let\casestyle=\relax
    \vbox{\smallskip\offinterlineskip
      \halign{&\vsquare{##}\cr #1}}}}

\def\bigyoung#1{\vcenter{%
    \squaresize=38pt\thickness=0.8pt\let\casestyle=\displaystyle
    \vbox{\smallskip\offinterlineskip
      \halign{&\vsquare{##}\cr #1}}}}

\title{\LARGE\sf
  Noncommutative symmetric functions \\
  and quasi-symmetric functions\\
  with two and more parameters}

\author{
  Florent {\sc Hivert},
  Alain {\sc Lascoux} and
  Jean-Yves {\sc Thibon}}

\date{}

\begin{document}

\maketitle

\begin{abstract}
We define two-parameter families of noncommutative symmetric functions and
quasi-symmetric functions, which appear to be the proper analogues of the
Macdonald symmetric functions in these settings.
\end{abstract}

\section{Introduction}

During the past twenty years, the theory of symmetric functions
underwent remarkable developments, including the discovery of
Macdonald's symmetric functions, and of various kinds of orthogonal
polynomials in several variables.

It also appeared that the classical theory of symmetric functions had
to be supplemented with the so-called noncommutative symmetric
functions and quasi-symmetric functions. These notions, originally
introduced for combinatorial purposes \cite{Ges,NCSF1,NCSF2}, are now
known to describe the representations of Hecke algebras and quantum
groups of type $A$ in a degenerate case, which is only apparent in a
certain parametrization (in the usual convention, the quantum
parameter $q$ stands for the point $[q:q^{-1}]$ of a complex
projective line, and the relevant degeneracy occurs at the origin
$[0:1]$ and at the point at infinity $[1:0]$)
\cite{NCSF3,NCSF4,NCSF5,NCSF6}.

Quasi-symmetric functions, whose algebra is denoted by $QSym$, and
noncommutative symmetric functions, denoted by $\Sym$, are two graded
Hopf algebras in natural duality. It has been found that almost all
interesting objects of the classical theory find an analogue on one
side or the other. There exist ``Schur functions'' in both set-ups.
Elementary and complete functions fall on the noncommutative side,
while monomial functions have only a quasi-symmetric analogue. Noncommutative
power sums are  not unique and turned out to be most
interesting from a combinatorial point of view. There are analogues of
the Littlewood-Richardson rule, of the Robinson-Schensted-Knuth
correspondence, of the internal product, of the Frobenius
characteristic map, of the Weyl and Demazure character formulae, and
so on.

It is therefore quite natural to push the analogy further, and to look
for analogues of the most interesting objects of the ``modern'' theory
of symmetric functions, i.e., Hall-Littlewood and Macdonald functions.

Although we have (at the time of writing) no analogue of the Hall
algebra to motivate the introduction of quasi-symmetric and
noncommutative Hall-Littlewood functions, recent developments
connecting these objects to Hecke algebras and Kazhdan-Lusztig
polynomials allows one to look for a definition involving a a
quasi-symmetrizing action of the Hecke algebra. This has been achieved
in \cite{Hiv}. The resulting quasi-symmetric Hall-Littlewood
functions, and their noncommutative dual basis have been shown to
share many properties with their classical analogues, which makes more
than plausible the pertinence of their definition.

In this note, we propose a combinatorial definition of noncommutative
and quasi-symmet\-ric analogues of Macdonald's symmetric functions. As
we shall see, the naive analogues of the
definitions by triangularity properties do not work (although
such a definition is in fact possible),
but the constraints that we should recover Hall-Littlewood functions
at $q=0$, and that the four-fold symmetry of the classical
$(q,t)$-Kostka matrix has to be retained (because of a conjectural
representation theoretical interpretation) are sufficient to suggest a
general pattern. Once the noncommutative analogues have been found,
the quasi-symmetric ones can be defined by duality.

The noncommutative theory is not expected to yield much information
about the classical Macdonald functions. What is expected, but
not proved yet, is that both represent different  projections of some
higher level object, which remains to be discovered. To support
this hypothesis, we derive a few properties of the noncommutative
and quasi-symmetric analogues which are direct analogues of known
properties of the classical Macdonald functions.
 
Details, proofs and other results will appear in a subsequent paper.
Our notation is as in \cite{NCSF1} and \cite{Mcd}

\footnotesize
{\it Acknowledgements} This research has been carried out at
the Isaac Newton Institute for Mathematical Sciences, during the
program {\it Symmetric functions and Macdonald polynomials}, whose
support is gratefully acknowledged. 
\normalsize

\newpage
\section{Noncommutative analogues of \\ the Macdonald functions}

We shall start with the definition of noncommutative analogues of the
symmetric functions
\begin{equation}
\tilde H_\mu(X;q,t)=\sum_\lambda \tilde K_{\lambda\mu}(q,t) s_\lambda(X)
=t^{n(\mu)}J_\mu\left( {X\over 1-t^{-1}};q,t^{-1}\right)\,.
\end{equation}
It has been conjectured by Garsia and Haiman \cite{GH}, and recently
proved by Haiman \cite{Hai} that these functions were the bigraded
Frobenius characteristics of certain realizations of the regular
representations of the symmetric group.

If we want to define noncommutative analogues $\tHH_J(A;q,t)$, to be 
called noncommutative Macdonald functions,
labelled by compositions $J$, as
\begin{equation}\label{eq:HtildaJ}
\tHH_J(A;q,t)= \sum_I \tilde k_{IJ}(q,t) R_I(A)
\end{equation}
where $R_I$ are the noncommutative ribbon Schur functions, it is
natural to invoke the representation theoretical interpretation of the
ribbons. It is known that they are the characteristics of the
indecomposable projective modules of the $0$-Hecke algebra $H_n(0)$,
each of them occuring with multiplicity one in the decomposition of
the regular representation. Hence, if we expect (\ref{eq:HtildaJ}) to
describe the bigraded characteristic of a regular representation of
$H_n(0)$, the $k_{IJ}(q,t)$ have to be monomials $q^it^j$. This will
be our first requirement.

Our second requirement is that 
\begin{equation}
\tHH_J(A;0,t)=t^{{l(J)\choose 2}} H_J(A;t^{-1})
\end{equation}
where in the right-hand side, $H_J$ is the noncommutative
Hall-Littlewood function of \cite{Hiv}.

Finally, the structure of the projective modules of $H_n(0)$ leads us
to expect that the $(q,t)$-Kostka monomials should possess the
symmetries
\begin{eqnarray}
\tilde k_{I \bar J^\sim}(q,t) &=& \tilde k_{IJ}(t,q)\,,\\
\tilde k_{IJ}(q,t)\tilde k_{\bar I^\sim J}(q,t)&=& q^{{n+1-l(J)\choose 2}}t^{{l(I)\choose 2}}\,,
\end{eqnarray}
and that $\tilde k_{(n),J}(q,t)$ is always equal to 1.

These requirements are sufficient to determine the first matrices,
whose transposes (denoted by $K_n$) are reproduced below.

$$
K_2=\begin{pmatrix}2 & 1 & q\cr 11 & 1 & t\end{pmatrix}
$$

$$
K_3=\begin{pmatrix}3& 1& q^{2}& q& q^{3}\cr {21}& 1& t& q& tq\cr 
    {12}& 1& q& t& tq\cr {111}& 1& t^{2}& t& t^{3}\end{pmatrix}
$$

$$
K_4=\begin{pmatrix}
  4& 1& q^{3}& q^{2}& q^{5}& q& q^{4}& q^{3}& q^{6}\cr 
{31}& 1& t& q^{2}& tq^{2}& q& tq& q^{3}& tq^{3}\cr {22}& 
1& q^{2}& t& tq^{2}& q& q^{3}& tq& tq^{3}\cr {211}& 1& t^{
2}& t& t^{3}& q& t^{2}q& tq& t^{3}q\cr {13}& 1& q^{2}& q& 
q^{3}& t& tq^{2}& tq& tq^{3}\cr {121}& 1& t^{2}& q& t^{2}
q& t& t^{3}& tq& t^{3}q\cr {112}& 1& q& t^{2}& t^{2}q& t& 
tq& t^{3}& t^{3}q\cr {1111}& 1& t^{3}& t^{2}& t^{5}& t& t
^{4}& t^{3}& t^{6}\end{pmatrix}
$$

\bigskip
These matrices have an apparent $2\times 2$-block structure. Actually,
$K_{n-1}$ is the submatrix of $K_n$ formed by even rows 
and odd columns. The rest of the matrix is determined by
its second column, corresponding to the composition $I=(n-1,1)$.
The monomials in this column are completely determined by our requirements
and it is found that $\tilde k_{(n-1,1),J}$ is equal to $q^{n-l(J)}$
if $J$ is an odd row, and to $t^{l(J)-1}$ if $J$ labels an even row.
Then, each pair of consecutive rows $(1,2),(3,4),\dots $ is labelled
by compositions $(L',L'')$, where $L$ is a composition of $n-1$,
$L'$ is obtained from $L$ by incrementing its last part of 1 and
$L''$ by adding at the end of $L$ a part 1, e.g., $(21)'=(22)$
and $(21)''=(211)$.
Let 
$$
A=\begin{pmatrix}1 & q^{l(L)} \cr 
           1 & t^{n-1-l(L)}\end{pmatrix}
$$
be the block in the firsts two columns. Then, the block
in the columns $(M',M'')$ is $\tilde k_{L,M}(q,t)A$.
This description will be taken as the definition of $K_n$
for general $n$.

The entries of the matrix $K_n$ can be directly described
in terms of the geometry of composition diagrams .

Given a composition $I$ of $n$, let $\Des(I)\subseteq
\{1,\ldots,n-1\}$ be the descent set of $I$.  For $1\leq k\leq n-1$,
put
$$d(I,k) = \# \{ k'<k,\, k'\in \Des(I) \}$$
and 
$$ v(I,k) =
\begin{cases}
  t^{1+d(I,k)} & \text{if $k     \in \Des(I)$}\,, \\
  q^{k-d(I,k)} & \text{if $k \not\in \Des(I)$}\,.
\end{cases}
$$
Then,
\begin{equation}\label{KIJ}
K_n(I,J) = \prod_{k\in \Des(J)} v(I,k)  \,.
\end{equation}

For $I=(1,2,2,4,1,1,3)$ and $J=(2,3,1,2,3,1,2)$ for example, one has
$$
\Des(I)=\{ 1,3,5,9,10,11 \}, \Des(J)= \{ 2,5,6,8,11,12\}
$$
and thus
$$ 
K_{14}(I,J)= q^{2-1}\, t^{1+2}\, q^{6-3}\, q^{8-3}\, t^{1+5}\, q^{12-6}
= q^{15}\, t^9  \ .
$$

The determinant of the matrix $K_n(q,t)$ can be computed, and the
result is very similar to what is obtained in the classical case:
\begin{equation}
\det K_n(q,t) = \prod_{m=1}^{n-1}\prod_{k=1}^m\left( t^{m+1-k}-q^k\right)^{2^{n-1-m}{m-1\choose k-1}}\,.
\end{equation}

The condition that the $\tHH$ functions should reduce to
Hall-Littlewood functions at $q=0$ imply the specialization
\begin{equation}
\tHH_J(A;0,1)=S^J(A)
\end{equation}
and the symmetries of the matrix imply that
\begin{equation}
\tHH_J(A;1,0)=S^{\bar J^\sim}\,.
\end{equation}
Also, setting $t=q^{-1}$ and clearing the denominators, we have
\begin{equation}
q^{l(J)-1\choose 2} \tHH_J(A;q,q^{-1}) \equiv R_J(A)\ ({\rm mod}\ q{\cal L})
\end{equation}
where ${\cal L}$ is the $\Z[q]$-lattice spanned by the ribbons
$R_I(A)$ in $\Sym$.

The expansions of the $\tHH$ functions in some other bases of $\Sym$
can also be given in closed form. Given a composition $J$ of $n$, we
associate to the $k$-th box (where $1\le k\le n-1$) of the diagram of
$J$ the polynomials
\begin{equation}
e(J,k)=
\begin{cases}
  t^i-1 & \text{if $k     \in \Des(J)$}\,, \\
  q^j-1 & \text{if $k \not\in \Des(J)$}\,.
\end{cases}
\end{equation}
if this box is located in row $i$ and column $j$ of the diagram. Then,
the expansion of $\tHH_J$ on the basis $\Lambda^I$ (products of
elementary functions) is given by
\begin{equation}
\tHH_I=\sum_I \prod_{k\not\in\Des(I)} e(I,k) \Lambda^I \,.
\end{equation}
For example, the filling of the diagram of $J=(2,2)$ is
$$
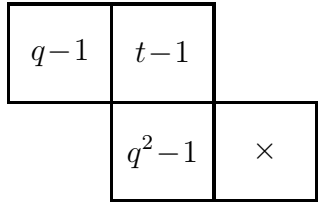
\bigyoung{ q\!-\!1 & t\!-\!1 \cr \blank &  q^2\!-\!1 & \times\cr}
$$
so that
\begin{eqnarray*}
\tHH_{22}=& \Lambda^{1111}+(q-1)\Lambda^{211}+(t-1)\Lambda^{121}
             +(q^2-1)\Lambda^{112} \\
& + (q-1)(t-1)\Lambda^{31}+(q-1)(q^2-1)\Lambda^{22}+(t-1)(q^2-1)\Lambda^{13}\\
&
  +(q-1)(t-1)(q^2-1)\Lambda^4\,.
\end{eqnarray*}

The expansion on the basis $(S^I)$ is, up to powers of $q$ and $t$,
given by the same formula.

\section{Quasi-symmetric analogues \\ of the Macdonald functions}

Quasi-symmetric Macdonald functions can now be defined by duality. The dual
basis of $(\tHH_J)$ in $QSym$ will be denoted by $(\tGG_I)$. We have
\begin{equation}
\tGG_I(X;q,t)= \sum_J \tilde g_{IJ}(q,t) F_J(X)
\end{equation}
where the coefficients are given by the transposed inverse
of the  Kostka matrix: $(\tilde g_{IJ})=(\tilde k_{IJ})^{-1}$. Remarkably,
there is for each $I$ a polynomial $D_I(q,t)$ such that
\begin{equation}
\tilde g_{IJ}(q,t)= (-1)^{l(I)-l(J)} {t^{a(I,J)}q^{b(I,J)}\over D_I(q,t)}
\end{equation}
It will be convenient to get rid of the common denominator and to introduce
the polynomials
\begin{equation}
\tPP_I(X;q,t)=D_I(q,t)\tGG_I(X;q,t)\,.
\end{equation}
The denominators are also common to pairs of consecutive columns. In the notation
of the previous section, if $L$ is a composition of $n-1$,
\begin{equation}
D_{L'}(q,t)=D_{L''}(q,t)\,.
\end{equation}
If $I=L'$ or $I=L''$, the formula is
\begin{equation}
D_I(q,t)=\prod_{(i,j)\in {\rm Diagr\,}(L)} (t^i-q^j)\,,
\end{equation}
where ${\rm Diagr\,}(L)$ is the ribbon diagram of the composition $L$, the
cells being labelled from top to bottom and left to right, i.e., $i$ is
the row number and $j$ the column number, as in a matrix. For example,
$$
{\rm Diagr\,}(3,1,2)=
\bigyoung{t_1-q_1 & t_1-q_2 & t_1-q_3 \cr
       \blank  & \blank  & t_2-q_3 \cr
      \blank   & \blank  & t_3-q_3 & t_3-q_4\cr}
$$

The numerators can be described as follows.
Define 
$$ u(I,k) =
\begin{cases}
  q^{ k-d(I,k)} & \text{if $k     \in \Des(I)$}\,, \\
  t^{d(I,k)}    & \text{if $k \not\in \Des(I)$} \,.
\end{cases}
$$

Then,
$$ t^{a(I,J)} q^{b(I,J)} = \prod_{k\not \in \Des(J)} u(I,k)  $$

For example, with the same $I$ and $J$ as before, one has
$$ \{1,\ldots, 13\} \setminus \Des(J)= \{1,3,4,7,9,10,13 \}$$
and the numerator of $\tilde g_{IJ}$ is 
$$(-1)^{6-7}(q^{1-0 }\, q^{3-1 }\, t^{1+2}\, t^{1+3}\, q^{9-3 }\,
q^{10-4 }\, t^{1+6}) = -t^{14}\, q^{15} $$

For $n=3$, the coefficient of $F_J$ in $\tPP_I$ is in row $J$ and column $I$
of the matrix

$$\begin{pmatrix}3& t^{2}& -tq^{2}& -t^{2}q& q^{2}\cr {21}& -t& t
& q& -q\cr {12}& -t& q^{2}& t^{2}& -q\cr {111}& 1& -1& 
-1& 1\end{pmatrix}
$$
and the denominators are
$$D_{3}=D_{21}= (t-q)(t-q^2)   $$
$$D_{12}=D_{111}= (t-q)(t^2-q) $$

For $n=4$, the matrix is

$$\begin{pmatrix}4& t^{3}& -t^{2}q^{3}& -t^{3}q^{2}& tq^{4}& -t^{4}q
& t^{2}q^{3}& t^{3}q^{2}& -q^{3}\cr {31}& -t^{2}& t^{2}& 
tq^{2}& -tq^{2}& t^{2}q& -t^{2}q& -q^{2}& q^{2}\cr {22}& 
-t^{2}& tq^{3}& t^{3}& -tq^{2}& t^{2}q& -q^{3}& -t^{3}q& q^{
2}\cr {211}& t& -t& -t& t& -q& q& q& -q\cr {13}& -t^{2
}& tq^{3}& t^{2}q^{2}& -q^{4}& t^{4}& -t^{2}q^{2}& -t^{3}q& 
q^{2}\cr {121}& t& -t& -q^{2}& q^{2}& -t^{2}& t^{2}& q& -
q\cr {112}& t& -q^{3}& -t^{2}& q^{2}& -t^{2}& q^{2}& t^{3
}& -q\cr {1111}& -1& 1& 1& -1& 1& -1& -1& 1\end{pmatrix}
$$
and the denominators are
\begin{align*}
D_{4}  &=D_{31}   = (t-q)(t  -q^2)(t  -q^3) \\
D_{22} &=D_{211}  = (t-q)(t  -q^2)(t^2-q^2)  \\
D_{13} &=D_{121}  = (t-q)(t^2-q)  (t^2-q^2) \\
D_{112}&=D_{1111} = (t-q)(t^2-1)  (t^3-q)
\end{align*}

\section{Multiparameter versions and the multiplication rule}

It is possible to replace $q^i$ and $t^j$ by independent
indeterminates $q_i,t_j$ in the combinatorial rule (\ref{KIJ}), and it
turns out that most of the previous formulas remain valid for the new
functions. Moreover, the multiplication rule is more transparent on
the multiparameter version.

\bigskip
Let $Z= \{z_0=1, z_1, z_2,\ldots \}$ be a sequence of commuting
indeterminates. We set
\begin{equation}
  \K_n(A;Z) = \sum_{|I|=n}
       \left( \prod_{d\in \Des(I)} z_d \right)     R_I \,.
\end{equation}
Given a composition $J$ of $n$, let  
$$ \vt(J,k) = 
\begin{cases}
  t_{1+d(J,k)} & \text{if $k     \in \Des(J)$}\,, \\
  q_{k-d(J,k)} & \text{if $k \not\in \Des(J)$}\,.
\end{cases}
$$
and let 
$$Z(J)  = \{ z_0=1, z_1=\vt(J,1), z_2=\vt(J,2),\ldots, z_{n-1}=\vt(J,n-1)  \} \,.$$

We now define the \emph{multiparameter noncommutative Macdonald
  functions} by: 
$$ \ttH_J(A;Q,T) =   \K_n(A;Z(J))  \ .$$
\bigskip
Thus $\tHH_J$ is the image of $\ttH_J$ under $t_j\mapsto t^j$,
$q_j\mapsto q^j$, $j=1,2,\ldots,n-1$.  The multiparameter Kostka
monomials $\bkt_{IJ}(Q,T)$ are defined as the coefficients in the
expansion of $\ttH_J$ in terms of ribbons
\begin{equation}
\ttH_J(A;Q,T)=\sum_I \bkt_{IJ}(Q,T) R_I(A)\,.
\end{equation}

When $Q=0$, one obtains in particular multiparameter
Hall-Littlewood functions. In the commutative case,
such multiparameter functions have been defined in \cite{LLT},
in the case of rectangular partitions, and their representation
theoretical meaning has been explained in \cite{Cal}. 
The noncommutative multiparameter
Hall-Littlewood function corresponding to the column
composition has been introduced in \cite{NCSF2}, and is
used in \cite{LLT} to give a closed expression for
the commutative one.

For $n=3,4$, the transposed Kostka matrices are as follows:

$$
\K_3=
\left (\begin {array}{ccccc}
3& 1&{ q_2}&{ q_1}&{ q_1}\,{ q_2}
\\
21&1&{ t_1}&{ q_1}&{ q_1}\,{ t_1}
\\
12&1&{ q_1}&{ t_1}&{ q_1}\,{ t_1}
\\
111&1&{ t_2}&{ t_1}&{ t_1}\,{ t_2}\end {array}
\right )
$$

$$
\K_4=
\left (\begin {array}{ccccccccc} 
4 &1&{ q_3}&{ q_2}&{ q_2}\,{ q_3}
&{ q_1}&{ q_1}\,{ q_3}&{ q_1}\,{ q_2}&{ q_1}\,{ q_2}\,
{  q_3}\\
31&1&{ t_1}&{ q_2}&{ q_2}\,{ t_1}&
{  q_1}&{ q_1}\,{ t_1}&{ q_1}\,{ q_2}&{ q_1}\,{ q_2}\,
{  t_1}\\
22&1&{ q_2}&{ t_1}&{ q_2}\,{ t_1}&{
 q_1}&{ q_1}\,{ q_2}&{ q_1}\,{ t_1}&{ q_1}\,{ q_2}\,{
 t_1}\\
211&1&{ t_2}&{ t_1}&{ t_1}\,{ t_2}&{
 q_1}&{ q_1}\,{ t_2}&{ q_1}\,{ t_1}&{ q_1}\,{ t_1}\,{
 t_2}\\
13&1&{ q_2}&{ q_1}&{ q_1}\,{ q_2}&{
 t_1}&{ q_2}\,{ t_1}&{ q_1}\,{ t_1}&{ q_1}\,{ q_2}\,{
 t_1}\\
121&1&{ t_2}&{ q_1}&{ q_1}\,{ t_2}&{
 t_1}&{ t_1}\,{ t_2}&{ q_1}\,{ t_1}&{ q_1}\,{ t_1}\,{
 t_2}\\
112&1&{ q_1}&{ t_2}&{ q_1}\,{ t_2}&{
 t_1}&{ q_1}\,{ t_1}&{ t_1}\,{ t_2}&{ q_1}\,{ t_1}\,{
 t_2}\\
1111&1&{ t_3}&{ t_2}&{ t_2}\,{ t_3}&{
 t_1}&{ t_1}\,{ t_3}&{ t_1}\,{ t_2}&{ t_1}\,{ t_2}\,{
 t_3}
\end {array}\right )
$$

The factorisation of the determinant still holds:

\begin{equation}
\det\K_n(Q,T) = \prod_{m=1}^{n-1}
                \prod_{k=1}^m
                     (t_{m+1-k}-q_k)^{ 2^{n-1-m}{m-1\choose k-1} }\,.
\end{equation}

When one set of parameters is specialized to 1, e.g.
$t_1=t_2=\ldots=t_{n-1}=1$,
one has the factorization
\begin{equation}
\ttH_J(A;{\bf 1},Q)=\ttH_{j_1}(A;{\bf 1},Q_1)\ttH_{j_2}(A;{\bf 1},Q_2)
  \cdots \ttH_{j_r}(A;{\bf 1},Q_r)\,,
\end{equation}
where $Q_1=\{q_1,\ldots,q_{j_1-1}\}$, $Q_2=\{q_{j_1},\ldots,q_{j_1+j_2-1}\}$, and
so on.

Moreover, the products
$$\K_n(Q,T)\K_n(Q,{\bf 1})^{-1} \ \ {\rm and}\ \ \K_n(Q,T)\K_n({\bf
  1},T)^{-1} $$
are respectively lower and upper triangular matrices,
whose diagonal elements can be given explicitly. This can be seen as
the noncommutative multiparameter  version of the characterization of
Macdonald's polynomials $\tHH_\mu(X;q,t)$ based on the two
$\lambda$-ring transformations $f(X)\rightarrow f((1-q)X)$ and
$f(X)\rightarrow f((1-t)X)$ (\textit{cf.} \cite{Hai}). There are several possible
tranformations extending $X\rightarrow (1-q)X$ in the noncommutative
case. The one which is involved here is different from the
transformation appearing in \cite{NCSF1,NCSF2}, even when $q_i=q^i$ or
$t_i=t^i$.

Given three integers $n,r,c$, associate to every composition 
$J$ of $n+1$ the following rational function read off from
the diagram of $J$ (boxes are numeroted from $0$ to $n$ and $z_0=1$).

Place the first box of the diagram of $J$ in position $(r,c)$.
Fill this box with $1/(t_r-q_c)$, and box number $k$, for
$k=1,\ldots, n$, in row $i$ and column $j$,
with
\begin{equation*}
  \begin{cases}
    \displaystyle\frac{t_i - z_{k-1}}{t_i -q_j} &
    \text{if the box has no other box above itself in its column}\\
    \displaystyle\frac{z_{k-1} - q_j}{t_i -q_j} &
    \text{otherwise.}\\
  \end{cases}
\end{equation*}
For the last box ($k=n$) omit the denominator. Let $\phi(J,r,c)$ be
the product of all the entries contained in the boxes of $J$.

With this at hand, we can now state the multiplication rule for
the $\ttH$-functions.
Let $I$ be a composition, $n$ an integer. 
Let $(r,c)$ be the coordinates of the last box of the diagram of $I$.
Then  
$$ \ttH_I(A;Q,T)\, \K_n(A;Z) =  
  \sum_J \phi(J', r,c)\, \ttH_J(A;Q,T)   \,,$$
where the sum is over the $2^{n}$ compositions $J$
obtained by concatenating $I$ with a composition of $n$,
$J'$ being the restriction of $J$ to its last $n+1$ boxes.

For example,  the product $\ttH_{12}\,  \K_2(A; 1,z_1) $ decomposes
into the following  sums of $\ttH_J(A;Q,T)$ functions
(writing the coefficients in the boxes of the diagram of $J$)

\begin{gather*}
\bigyoung{ 
  \cr
  & \frac{1}{t_2-q_2}   & \frac{t_2-1}{t_2-q_3}  & \frac{t_2-z_1}{1}  \cr } 
\quad + \quad
\bigyoung{ 
  \cr
  & \frac{1}{t_2-q_2}   & \frac{t_2-1}{t_2-q_3} \cr
  \blank &\blank &\frac{z_1-q_3}{1}    \cr }
\\
\quad + \quad
\bigyoung{ 
  \cr
  & \frac{1}{t_2-q_2}  \cr
  \blank & \frac{1-q_2}{t_3-q_2} & \frac{t_3-z_1}{1} \cr}
\quad + \quad
\bigyoung{ 
  \cr
  & \frac{1}{t_2-q_2}  \cr
  \blank & \frac{1-q_2}{t_3-q_2}\cr
  \blank & \frac{z_1-q_2}{1}    \cr }
\end{gather*}
which reads
\begin{align*}
 \ttH_{12}\,  \K_2(A; 1,z_1)\  = 
\     {({t_2}-1)({t_2}-{z_1})\over ({t_2}-{q_2})({t_2}-{q_3})}\,\ttH_{14} 
\ +\ &{({t_2}-1)({z_1}-{q_3})\over ({t_2}-{q_2})({t_2}-{q_3})}\,\ttH_{131} \\
\ +\  {(1-{q_2})({t_3}-{z_1})\over ({t_2}-{q_2})({t_3}-{q_2})}\,\ttH_{122} 
\ +\ &{(1-{q_2})({z_1}-{q_2})\over ({t_2}-{q_2})({t_3}-{q_2})}\,\ttH_{1211}  
\end{align*}

\bigskip
{\footnotesize\noindent
\sc
Isaac Newton Institute for Mathematical Science\\
20 Clarkson road\\
Cambridge, CB3 0EH, U.K.\\
\\
\sf Permanent address:\\
\sc Institut Gaspard Monge, \\
Universit\'e de Marne-la-Vall\'ee,\\
77454 Marne-la-Vall\'ee cedex,\\
 FRANCE}

\end{document}